\documentclass[10pt,twoside]{article}

\usepackage[centertags]{amsmath}
\usepackage{amsfonts}
\usepackage{amssymb}
\usepackage{amsthm}
\usepackage{epsfig}
\usepackage{color}
\allowdisplaybreaks[4]
\definecolor{c20}{rgb}{0.,0.7,0.}
\definecolor{c30}{rgb}{0.,0.,1.}
\definecolor{c40}{rgb}{1,0.1,0.7}
\definecolor{c50}{rgb}{1,0,0}
\definecolor{c60}{rgb}{1,0.9,0.1}

\date{}

\begin{document}
\baselineskip 15pt \setcounter{page}{1}
\title{\bf \Large  The limit theorems on extreme order statistics and partial sums of i.i.d. random variables
\thanks{Research supported by Innovation of Jiaxing City: a program to support the talented persons
and Project of new economy research center of Jiaxing City (No. WYZB202254).}}
\author{{\small Gaoyu Li$^{1,2}$ \ \ Zhongquan Tan$^{2}$\footnote{Corresponding author. E-mail address: tzq728@163.com }}\\
\\
{\small\it 1. Department of Mathematic, Zhejiang Normal University, Jinhua 321004, PR China }\\
{\small\it 2. College of Data Science, Jiaxing University, Jiaxing 314001, PR China}\\
}
 \maketitle
 \baselineskip 15pt

\begin{quote}
{\bf Abstract:}\ \  This paper proves several weak limit theorems for the joint version of extreme order statistics and partial sums of independently and identically distributed random variables. The results are also extended to almost sure limit version.



\end{quote}

\section{Introduction }
Since the maxima and partial sums of a sequence of random variables played very important role in theoretical and applied probability,
many authors studied their joint asymptotic properties.
Chow and Teugels (1978) first investigated asymptotic relation between the maxima and partial sums of a sequence of independently and identically distributed (i.i.d.) random variables $\{X_{n}, n\geq1\}$.  Denote by  $M_{n}$ and $S_{n}$ the maximum and partial sum of $\{X_{1},X_{2},\cdots, X_{n}\}$, respectively.
Suppose that there exist constants $a_{n}>0,b_{n}\in \mathbb{R}, n\ge1$ and a non-degenerate distribution $G(x)$ such that
\begin{eqnarray*}
\label{eq1.1}
\lim_{n\to\infty}P(a_{n}(M_{n}-b_{n})\le x)=G(x),
\end{eqnarray*}
where $G$ must be one of the following three types of extreme value distributions:
\begin{eqnarray*}
\mbox{Gumbel:}\ \ \Lambda(x)&=&\exp\left(-e^{-x}\right),\ \ -\infty<x<+\infty;\\
\mbox{Frechet:}\ \ \Phi_{\alpha}(x)&=&\left\{\begin{array}{l}0, \quad\quad\quad\quad\quad\ \ x\leq 0, \\ \exp \left(-x^{-\alpha}\right), \quad x>0,\end{array}\right.\alpha\geq 0;\\
\mbox{Weibull:}\ \ \Psi_{\alpha}(x)&=&\left\{\begin{array}{l}\exp \left(-(-x)^{\alpha}\right), \quad x\leq0, \\ 1,\ \ \quad\quad\quad\quad\quad\quad x >0, \end{array}\right.\alpha>0.
\end{eqnarray*}
If $\{X_{n}, n\geq1\}$ has finite variance, Chow and Teugels (1978) showed that
\begin{equation}
\label{eq1.2}
\lim_{n\rightarrow\infty}P\left(a_{n}(M_{n}-b_{n})\leq x, \frac{S_{n}}{\sqrt{n}}\leq y\right)=G(x)\Phi(y),
\end{equation}
which indicates the maximum and the partial sum are asymptotically independent.
If $\{X_{n}, n\geq1\}$ has infinite variance, Chow and Teugels (1978) showed that the maximum and the partial sum are asymptotically dependent.

The result of (\ref{eq1.2}) has many extensions. A natural extension is to consider the dependent case. Anderson and Turkman (1991a,1991b) extended Chow and Teugels's result to strong mixing case under a technical condition
and then Hsing (1995) completed their results by deleting the technical condition.
Ho and Hsing (1996), Ho and McCormick (1999), McCormick and Qi (2000) and Peng and Nadarajah (2002) dealt with this problem for some dependent Gaussian cases and
James et al. (2007) considered the multivariate stationary Gaussian cases.

By studying the joint behavior of the point processes of exceedances and the partial sums for Gaussian processes, Hu et al. (2009) derived the
asymptotic distribution of the extreme order statistics and partial sums.
Let $\{X_{n}, n\geq1\}$  be a stationary Gaussian sequence with covariance functions $r(n)=E(X_{1}X_{n+1})$.
Denote by  $M_{n}^{(k)}$ the $k$-th maximum of $\{X_{1},X_{2},\cdots, X_{n}\}$.
Suppose that the well known Berman's condition was satisfied, i.e.,
$$r(n)\log n\rightarrow 0\ \ \mbox{as}\ \ n\rightarrow \infty,$$
then for any  $x, y\in \mathbb{R}$,
\begin{equation}
\label{eq1.3}
\lim_{n\rightarrow\infty}P\left(c_{n}(M_{n}^{(k)}-d_{n})\leq x, \frac{S_{n}}{\sqrt{Var(S_{n})}}\leq y\right)=\exp(-e^{-x})\sum_{i=0}^{k-1}\frac{(e^{-x})^{i}}{i!}\Phi(y),
\end{equation}
where $\Phi (\cdot)$ stands for the standard normal distribution function and  $c_{n}>0, d_{n}\in \mathbb{R}, n\ge
1$ are real sequences.
For more studies in this direction, we refer to Peng (1998), Tan and Peng (2011a, b), Peng et al. (2012) and Tan and Yang (2015).

Another interesting extension is to consider the almost sure limit theorem for the maxima and partial sums.
Based on the result of Chow and Teugels (1978), Peng et al. (2009) studied the almost sure limit theorem for the maxima and partial sums
for i.i.d. random variables and showed that
\begin{equation}
\label{eq1.4}
 \lim_{N\rightarrow \infty}\frac{1}{\log
N}\sum^{N}_{n=1}\frac{1}{n}I\left(a_{n}(M_{n}-b_{n})\le
x,\frac{S_{n}}{\sqrt{n}}\leq y\right)=G(x)\Phi(y)\ \ a.s.
\end{equation}
for any $x,y\in \mathbb{R}$, where $I(\cdot)$ denotes the indicator function.
For more studies on almost sure limit theorem for the maxima and partial sums for dependent cases, we refer to Dudzi\'{n}ski (2003, 2008),  Zang et al. (2010), Zhao et al. (2011), Tan and Wang (2011), Zang (2012),  Wu (2015) and Wu and Jiang (2016).

Most of the above mentioned studies focus on the Gaussian cases, especially regarding the joint version between extreme order statistics and partial sums.
We even don't know the asymptotic relation between extreme order statistics and partial sums  for i.i.d cases.
In this paper, we continue to study the joint asymptotic properties of extreme order statistics and partial
sums for i.i.d. random sequence and have two goals. The first one is to derive the joint asymptotic distribution of
$k$-th maxima and partial sums, which will extend the classic result (\ref{eq1.2}). The second one is to prove the almost sure limit theorem
for the joint version of $k$-th maxima and partial sums, which will extend the result (\ref{eq1.4}).
Some related limit results on the joint asymptotic of the $k$ largest maxima and partial sums are also given.

\section{The weak limit theorems}
For any random sequence $\{X_{n}, n\geq1\}$,  define the point process of
exceedances of level $u_{n}$ formed by $X_{i}, 1\leq i\leq n$ as
\begin{equation}
\label{eq2.1}
N_{n}^{X}(B)=\sum_{i=1}^{n}I(X_{i}>u_{n}, i/n \in B)
\end{equation}
for any Borel set $B$ on $(0,1]$.

In this paper, let $M_{n}^{(k)}$ be the $k$-th maximum of $\{X_{1},X_{2},\cdots, X_{n}\}$ and $S_{n}=\sum_{i=1}^{n}X_{i}$.
Let $u_{n}=u_{n}(x)=a_{n}^{-1}x+b_{n}$, where $a_{n}>0,b_{n}\in \mathbb{R}, n\ge1$ will be chosen in Theorem 2.1.

Now we state our main results.

\textbf{Theorem 2.1}. {\sl Let $\{X_{n}, n\geq1\}$  be a sequence of i.i.d. random variables with nondegenerate common distribution function $F$ and $E(X_{1})=0$, $E(X_{1}^{2})=1$. Suppose that there exist constants $a_{n}>0,b_{n}\in \mathbb{R}, n\ge1$  such that
 \begin{equation}
\label{eq2.1}
\lim_{n\to\infty}P(a_{n}(M_{n}-b_{n})\le x)=G(x).
\end{equation}
 holds, where $G$ is one of the extreme value distribution $\Lambda$, $\Phi_{\alpha}$ for some $\alpha>2$ and $\Psi_{\alpha}$ for some $\alpha>0$.
Then $N_{n}^{X}\stackrel{d}\rightarrow N$ and
$N_{n}^{X}$ and $\frac{S_{n}}{\sqrt{n}}$ are asymptotically
independent, where $N$ is a Poisson point process on $(0,1]$ with intensity $\log G^{-1}(x)$.}

\textbf{Remark 2.2}. {\sl  This paper only deals with the case that $F$ has finite variance. However, our method also does not hold for the case
that $F$ is in the domain attraction of $\Phi_{\alpha}$  with $\alpha=2$.
}

\textbf{Theorem  2.3}. {\sl Under the conditions of Theorem 2.1, we have for any $x,y\in \mathbb{R}$
\begin{equation}
\label{eq2.2}
\lim_{n\rightarrow\infty}P\left(a_{n}(M_{n}^{(k)}-b_{n})\leq x,\frac{ S_{n}}{\sqrt{n}}\leq y\right)=G(x)\sum_{i=0}^{k-1}\frac{(-\ln G(x))^{i}}{i!}\Phi(y).
\end{equation}
}

The above results can be extended to the point process of exceedances of multiple levels.
Let $u_{n}^{(1)}\geq u_{n}^{(2)}\geq,\ldots,u_{n}^{(s)}$ be $s$ levels, where $u_{n}^{(j)}=u_{n}(x_{j})=a_{n}^{-1}x_{j}+b_{n}$, $x_{j}\in \mathbb{R}$, $j=1,2,\ldots,s$.  Define the point process $\widetilde{N}_{n}$ of exceedances of levels $u_{n}^{(1)},u_{n}^{(2)},\ldots,u_{n}^{(s)}$ as
\begin{equation}
\label{eqt2.1}
\widetilde{N}_{n}^{X}(B)=\sum_{j=1}^{s}N_{n}^{(j)}(B)=\sum_{j=1}^{s}\sum_{i=1}^{n}I(X_{i}>u_{n}^{j}, (i/n,j) \in B)
\end{equation}
for any Borel set $B$ on $(0,1]\times \mathbb{R}$.

Construct a point process $\widetilde{N}$ on $(0,1]\times \mathbb{R}$ such that $\widetilde{N}(\cdot)=\sum_{j=1}^{s}N^{(j)}(\cdot)$, where
for each $j=1,2,\ldots, s-1$, $N^{(j)}$ is a Poisson process independent thinning of the Poisson process $N^{(j+1)}$ with deleting probability
$1-\frac{\log G^{-1}(x_{j})}{\log G^{-1}(x_{j+1})}$, the initial Poisson process $N^{(s)}$ has intensity $\log G^{-1}(x_{s})$.
 It is known that $N^{(j)}, j=1,2,\ldots, s$ are Poisson processes with intensity $\log G^{-1}(x_{j})$ respectively, and are independent on disjoint
intervals on the plane  (cf. Leadbetter et al. 1983, Section 5.5).

\textbf{Theorem 2.4}. {\sl Let $\{X_{n}, n\geq1\}$  be a sequence of standard i.i.d. random variables with nondegenerate common distribution function $F$ and $E(X_{1})=0$, $E(X_{1}^{2})=1$. Suppose that there exist constants $a_{n}>0,b_{n}\in \mathbb{R}, n\ge1$ and a non-degenerate distribution $G(x)$ such that
(\ref{eq2.1}) holds, where $G$ is one of the extreme value distribution $\Lambda$, $\Phi_{\alpha}$ for some $\alpha>2$ and $\Psi_{\alpha}$ for some $\alpha>0$.
Then $\widetilde{N}_{n}^{X}\stackrel{d}\rightarrow \widetilde{N}$ and
$\widetilde{N}_{n}^{X}$ and $\frac{S_{n}}{\sqrt{n}}$ are asymptotically
independent.}

It is easy to see from Theorem 2.4 that the $k$ largest maxima and the partial sum are asymptotically independent.

\textbf{Theorem 2.5}. {\sl Under the conditions of Theorem 2.4, we have for any fixed $k\geq 1$ and any $x_{1},\ldots,x_{k},y\in \mathbb{R}$
\begin{eqnarray}
\label{eq2.5}
&&\lim_{n\rightarrow\infty}P\left(a_{n}(M_{n}^{(1)}-b_{n})\leq x_{1},\ldots, a_{n}(M_{n}^{(k)}-b_{n})\leq x_{k},\frac{ S_{n}}{\sqrt{n}}\leq y\right)\nonumber\\
&&=H(x_{1},x_{2},\ldots,x_{k})\Phi(y),
\end{eqnarray}
where $H(x_{1},x_{2},\ldots,x_{k})$ is the joint limit distribution of $(a_{n}(M_{n}^{(1)}-b_{n}),\ldots,a_{n}(M_{n}^{(k)}-b_{n}))$.
}

The distribution function $H(x_{1},x_{2},\ldots,x_{k})$ has a very complicated representation. The following corollary states a two dimensional case.

\textbf{Corollary 2.6}. {\sl Under the conditions of Theorem 2.4, we have for any fixed $l>k\geq 1$ and any $x_{1},x_{2},y\in \mathbb{R}$
\begin{eqnarray}
\label{eq2.4}
&&\lim_{n\rightarrow\infty}P\left(a_{n}(M_{n}^{(k)}-b_{n})\leq x_{1}, a_{n}(M_{n}^{(l)}-b_{n})\leq x_{2},\frac{ S_{n}}{\sqrt{n}}\leq y\right)\nonumber\\
&&=\sum_{i=0}^{k-1}\sum_{j=i}^{l-1}G(x_{2})\frac{(\log G^{-1}(x_{1}))^{i}(\log G^{-1}(x_{2})-\log G^{-1}(x_{1}))^{j-i}}{i!(j-i)!}\Phi(y).
\end{eqnarray}
}

We need the following lemmas to prove Theorems 2.1 and 2.3-2.5.

\textbf{Lemma 2.7}. {\sl Let $\{X_{n}, n\geq1\}$  be a sequence of i.i.d. random variables with nondegenerate common distribution function $F$ and $E(X_{1})=0$, $E(X_{1}^{2})=1$. Let $\{Y_{n}, n\geq1\}$  be an independent copy of $\{X_{n}, n\geq1\}$. If $(\ref{eq2.1})$ holds, then for any $0<a<b\leq 1$
$$
\sum_{i=[na]+1}^{[nb]}\left|P\left(X_{i}> u_{n},
\frac{S_{n}}{\sqrt{n}}> y\right)- P\left(Y_{i}> u_{n},
\frac{S_{n}}{\sqrt{n}}> y\right)\right|\rightarrow 0,
$$
as $n\rightarrow \infty$.
}

\textbf{Proof:}  Note that $G$ defined in $(\ref{eq2.1})$ is either a Gumbel distribution $\Lambda(x)$ or a Frechet distribution $\Phi_{\alpha}(x)$ with $\alpha> 2$ or a Weibull-type distribution $\Psi_{\alpha}$ with $\alpha>0$.
Define $x_{F}=\sup\{x: F(x)< 1\}$ and $S_{ni}=S_{n}-X_{i}$, $i=1,2,\ldots,n$. Integration by parts shows that
\begin{eqnarray*}
P\left(X_{i}> u_{n},
\frac{S_{n}}{\sqrt{n}}> y\right)&=&\int_{u_{n}}^{x_{F}}P\left(
\frac{S_{ni}+x}{\sqrt{n}}> y\right)dF(x)\\
&=&(1-F(u_{n}))P\left(\frac{S_{ni}+u_{n}}{\sqrt{n}}> y\right)\\
&+&(1-F(u_{n}))\int_{u_{n}}^{x_{F}}\frac{1-F(x)}{1-F(u_{n})}dP\left(
\frac{S_{ni}+x}{\sqrt{n}}> y\right).
\end{eqnarray*}
Thus, the sum in Lemma 2.7 is bounded above by
\begin{eqnarray*}
&&\sum_{i=[na]+1}^{[nb]}\left|(1-F(u_{n}))P\left(\frac{S_{ni}+u_{n}}{\sqrt{n}}> y\right)-(1-F(u_{n}))P\left(\frac{S_{n}}{\sqrt{n}}> y\right)\right|\\
&&\ \ \ +\sum_{i=[na]+1}^{[nb]}(1-F(u_{n}))\int_{u_{n}}^{x_{F}}\frac{1-F(x)}{1-F(u_{n})}dP\left(
\frac{S_{ni}+x}{\sqrt{n}}> y\right)\\
&&=([nb]-[na])(1-F(u_{n}))\left|P\left(\frac{S_{n1}+u_{n}}{\sqrt{n}}> y\right)-P\left(\frac{S_{n}}{\sqrt{n}}> y\right)\right|\\
&&\ \ \ +([nb]-[na])(1-F(u_{n}))\int_{u_{n}}^{x_{F}}\frac{1-F(x)}{1-F(u_{n})}dP\left(
\frac{S_{n1}+x}{\sqrt{n}}> y\right)\\
&&=:([nb]-[na])(1-F(u_{n}))[T_{n1}+T_{n2}],
\end{eqnarray*}
where we used the stationarity of $\{X_{n}, n\geq1\}$ in the second step. Taking into account (\ref{eq2.1}), we get
\begin{eqnarray}
\label{Fact}
n(1-F(u_{n}))\rightarrow
\log G^{-1}(x)
\end{eqnarray}
as $n\rightarrow\infty$. Thus, to prove the lemma, it suffices to prove $T_{ni}\rightarrow0$, $i=1,2$, as $n\rightarrow\infty$.
If $G=\Lambda$  then $u_{n}$ is slowly varying  and so $u_{n}/\sqrt{n}\rightarrow0$, as $n\rightarrow\infty$ (see e.g., Anderson and Turkman (1991a)).
If $G=\Phi_{\alpha}$ with $\alpha>2$ then $u_{n}$ is regularly varying with index $1/\alpha<1/2$ and so again $u_{n}/\sqrt{n}\rightarrow0$, as $n\rightarrow\infty$ (see e.g., Anderson and Turkman (1991a)).
If $G=\Psi_{\alpha}$  with $\alpha>0$ then $u_{n}\leq x_{F}<\infty$ for all $n\geq 1$ and obviously $u_{n}/\sqrt{n}\rightarrow0$, as $n\rightarrow\infty$.
 Hence, by the central limit theorem for i.i.d random variables, $T_{n1}\rightarrow0$,  as $n\rightarrow\infty$.
It follows from the above arguments that we can also choose a constant $\varepsilon>0$ such that
$(u_{n})^{1+\varepsilon}/\sqrt{n}\rightarrow0$, as $n\rightarrow\infty$.
For the term $T_{n2}$, we will split it into two cases: $x_{F}=\infty$ and $x_{F}<\infty$.
For the case $x_{F}=\infty$, $G$ must be either $\Lambda$ or $\Phi_{\alpha}$.
Note that
\begin{eqnarray*}
T_{n2}&=&\int_{1}^{\infty}u_{n}\frac{1-F(tu_{n})}{1-F(u_{n})}dP\left(
\frac{S_{n1}+tu_{n}}{\sqrt{n}}> y\right).
\end{eqnarray*}
For $X$ in the domain of attraction of $\Phi_{\alpha}$, $\frac{1-F(tu_{n})}{1-F(u_{n})}\sim t^{-\alpha}$ as $u_{n}\rightarrow\infty$.
For $X$ in the domain of attraction of $\Lambda$, $\frac{1-F(tu_{n})}{1-F(u_{n})}\leq t^{-\alpha}$ for large enough $u_{n}$ whatever $\alpha>0$.
We thus have for $\alpha>2$
\begin{eqnarray*}
T_{n2}&\leq&\int_{1}^{\infty}u_{n}t^{-\alpha}dP\left(
\frac{S_{n1}+tu_{n}}{\sqrt{n}}> y\right)\\
&=&\int_{u_{n}}^{\infty}(\frac{x}{u_{n}})^{-\alpha}dP\left(
\frac{S_{n1}+x}{\sqrt{n}}> y\right)\\
&=&\int_{u_{n}}^{u_{n}^{1+\varepsilon}}(\frac{x}{u_{n}})^{-\alpha}dP\left(
\frac{S_{n1}+x}{\sqrt{n}}> y\right)+\int_{u_{n}^{1+\varepsilon}}^{\infty}(\frac{x}{u_{n}})^{-\alpha}dP\left(
\frac{S_{n1}+x}{\sqrt{n}}> y\right)\\
&\leq&\int_{u_{n}}^{u_{n}^{1+\varepsilon}}1dP\left(
\frac{S_{n1}+x}{\sqrt{n}}> y\right)+\int_{u_{n}^{1+\varepsilon}}^{\infty}(u_{n})^{-\varepsilon\alpha}dP\left(
\frac{S_{n1}+x}{\sqrt{n}}> y\right)\\
&\leq&P\left(\frac{S_{n1}+u_{n}^{1+\varepsilon}}{\sqrt{n}}> y\right)-P\left(\frac{S_{n1}+u_{n}}{\sqrt{n}}> y\right)
+(u_{n})^{-\varepsilon\alpha},
\end{eqnarray*}
which tends to $0$, by central limit theorem for i.i.d. random variables again, since $(u_{n})^{1+\varepsilon}/\sqrt{n}\rightarrow0$ and $u_{n}\rightarrow\infty$, as $n\rightarrow\infty$.
For the case $x_{F}<\infty$, $G$ must be either $\Lambda$ or $\Psi_{\alpha}$.
We have
\begin{eqnarray*}
T_{n2}&=&
\int_{u_{n}}^{x_{F}}\frac{1-F(x)}{1-F(u_{n})}dP\left(
\frac{S_{n1}+x}{\sqrt{n}}> y\right)\\
&\leq& P\left(\frac{S_{n1}+x_{F}}{\sqrt{n}}> y\right)-P\left(\frac{S_{n1}+u_{n}}{\sqrt{n}}> y\right)\rightarrow0,
\end{eqnarray*}
as $n\rightarrow\infty$, by central limit theorem for i.i.d. random variables again, since $u_{n}\leq x_{F}<\infty$ for all $n\geq 1$. The proof of the lemma is complete.

\textbf{Lemma 2.8}. {\sl Let $\{X_{n}, n\geq1\}$  be a sequence of i.i.d. random variables with nondegenerate common distribution function $F$ and $E(X_{1})=0$, $E(X_{1}^{2})=1$.
Let $M_{n}^{X}(a,b]=\max(X_{i},i/n \in (a,b])$ for any $0<a<b\leq 1$.
 If $(\ref{eq2.1})$ holds, then for any $0<a_{1}<b_{1}\leq a_{2}< b_{2}\leq \cdots \leq a_{k}<b_{k}\leq1$
$$
P\left(\bigcap_{i=1}^{k}\{M_{n}^{X}(a_{i},b_{i}]\leq u_{n}\},
\frac{S_{n}}{\sqrt{n}}\leq y\right)\to
\prod^{k}_{i=1}[G(x)]^{(b_{i}-a_{i})}\Phi(y),
$$}
as $n\rightarrow \infty$.

\textbf{Proof:} Since $\{X_{n}, n\geq1\}$ is a sequence of i.i.d. random variables, by Theorem 5.4.5 of Leadbetter et al. (1983),
we have, as $n\rightarrow\infty$
$$P\left(\bigcap_{i=1}^{k}\{M_{n}^{X}(a_{i},b_{i}]\leq u_{n}\}\right)\to
\prod^{k}_{i=1}[G(x)]^{(b_{i}-a_{i})}.$$
The asymptotic independence between the maxima and the sum can be proved by the same way as the proof of Theorem 1 of Anderson and  Turkman (1991a).

\textbf{Lemma 2.9}. {\sl Let $\{Y_{n}, n\geq1\}$  be an independent copy of $\{X_{n}, n\geq1\}$.
For any fixed $y\in \mathbb{R}$, constructing a new sequence as
$$
Z_{i}=X_{i}I(\frac{S_{n}}{\sqrt{n}}\leq
y)+{Y_{i}}I(\frac{S_{n}}{\sqrt{n}}>y),\ \ 1\leq i\leq n
$$
and define
\begin{equation}
\label{eq2.9.1}
N_{n}^{Z}(B)=\sum_{i=1}^{n}I(Z_{i}>u_{n}, i/n \in B)
\end{equation}
for any Borel set $B$ on $(0,1]$.
Then  $N_{n}^{Z}\stackrel{d}\rightarrow N$, where $N$ is a Poisson point process on $(0,1]$  with intensity $\ln G^{-1}(x)$.
}

\textbf{Proof:}  By Theorem A.1 in Leadbetter et al. (1983), it is sufficient to show that\\
$$ (i)\ \  EN^{Z}_{n}((a,b])\rightarrow EN((a,b])=(b-a)\ln G^{-1}(x), 0<a<b\leq 1;$$
\begin{eqnarray*}
(ii)\ \ P\left(\bigcap_{i=1}^{k}\{N^{Z}_{n}((a_{i},b_{i}])=0\}\right)
&\rightarrow & P\left(\bigcap_{i=1}^{k}\{N((a_{i},b_{i}])=0\}\right)\\
&=& \exp\left(-\sum_{i=1}^{k}(b_{i}-a_{i})\log G^{-1}(x)\right)
\end{eqnarray*}
where $0<a_{1}<b_{1}\leq a_{2}< b_{2}\leq \cdots \leq
a_{k}<b_{k}\leq 1$.

For (i), we have
\begin{eqnarray*}
EN^{Z}_{n}((a,b])&=&E\sum_{i/n \in (a,b]}I(Z_{i}>u_{n})\\
                 &=&\sum_{i=[na]+1}^{[nb]}\left\{P\left(Z_{i}>u_{n},\frac{S_{n}}{\sqrt{n}}\leq y\right)+
                    P\left(Z_{i}>u_{n},\frac{S_{n}}{\sqrt{n}}> y\right)\right\} \\
                 &=&\sum_{i=[na]+1}^{[nb]}\left\{P\left(X_{i}>u_{n},\frac{S_{n}}{\sqrt{n}}\leq y\right)+
                    P\left(Y_{i}>u_{n},\frac{S_{n}}{\sqrt{n}}> y\right)\right\} \\
                 &=&\sum_{i=[na]+1}^{[nb]}\left\{P(X_{i}>u_{n})+ P\left(Y_{i}> u_{n},\frac{S_{n}}{\sqrt{n}}> y\right)-
                    P\left(X_{i}> u_{n},\frac{S_{n}}{\sqrt{n}}> y\right)\right\} \\
                 &=&([nb]-[na])(1-F(u_{n})) \\
                 &&+\sum_{i=[na]+1}^{[nb]}\left\{P\left(Y_{i}> u_{n},\frac{S_{n}}{\sqrt{n}}> y\right)-
                    P\left(X_{i}> u_{n},\frac{S_{n}}{\sqrt{n}}> y\right)\right\}.
\end{eqnarray*}
 Thus by Lemma 2.7 and the fact (\ref{Fact}), we have
$$EN^{Z}_{n}((a,b])\rightarrow (b-a)\ln G^{-1}(x)=EN((a,b]),$$
as $ n\rightarrow\infty.$

Next for (ii), it follows  that
\begin{eqnarray}
\label{Lem2.1}
&&P\left(\bigcap_{i=1}^{k}\{N^{Z}_{n}((a_{i},b_{i}])=0\}\right)\nonumber\\
&&=P\left(\bigcap_{i=1}^{k}\{M_{n}^{X}(a_{i},b_{i}]\leq u_{n}\},
\frac{S_{n}}{\sqrt{n}}\leq y\right)
+P\left(\bigcap_{i=1}^{k}\{M_{n}^{Y}(a_{i},b_{i}]\leq
u_{n}\}\right)P\left(\frac{S_{n}}{\sqrt{n}}> y\right).
\end{eqnarray}
By Lemma 2.8, we have
\begin{eqnarray}
\label{Lem2.2}
P\left(\bigcap_{i=1}^{k}\{M_{n}^{X}(a_{i},b_{i}]\leq u_{n}\},
\frac{S_{n}}{\sqrt{n}}\leq y\right)
  \rightarrow \prod^{k}_{i=1}\exp(-(b_{i}-a_{i})\log G^{-1}(x))\Phi(y),
  \end{eqnarray}
as $ n\rightarrow\infty.$ Recalling that $\{Y_{n}, n\geq1\}$ is a sequence of
 i.i.d. random variables and noting the fact (\ref{Fact}) again,
we get
\begin{eqnarray}
\label{Lem2.3}
P\left(\bigcap_{i=1}^{k}\{M_{n}^{Y}(a_{i},b_{i}]\leq
u_{n}\}\right)\rightarrow \prod^{k}_{i=1}\exp(-(b_{i}-a_{i})\log G^{-1}(x)),
  \end{eqnarray}
as $ n\rightarrow\infty.$ Since $\{X_{n}, n\geq1\}$ is also a sequence of i.i.d. random variables with $E(X_{1})=0$, $E(X_{1}^{2})=1$, we have
\begin{eqnarray}
\label{Lem2.4}
P\left(\frac{S_{n}}{\sqrt{n}}> y\right)\rightarrow 1-\Phi(y),
\end{eqnarray}
as $ n\rightarrow\infty.$
Thus, plugging  (\ref{Lem2.2}-\ref{Lem2.4}) into (\ref{Lem2.1}), we get
$$P\left(\bigcap_{i=1}^{k}\{N^{Z}_{n}((a_{i},b_{i}])=0\}\right)\rightarrow\exp\left(-\sum_{i=1}^{k}(b_{i}-a_{i})\log G^{-1}(x)\right).$$
This completes the proof of Lemma 2.9.

\textbf{Proof of Theorem 2.1:} Since $\{X_{n}, n\geq1\}$ is a sequence of
standard i.i.d. random variables, by Theorem 5.2.1 in Leadbetter et al. (1983), $N_{n}^{X}$ converges in
distribution to a Poisson process on (0,1] with intensity $\log G^{-1}(x)$. Next, we show the asymptotic independence between $N_{n}^{X}$
and $S_{n}$.

First noting that for disjoint Borel sets $B_{1}, B_{2},\cdots, B_{k}$ on $(0,1]$
with $m(\partial B_{i})=0,i=1,2,\cdots,k$, we have
\begin{eqnarray*}
&&P\left(\bigcap_{i=1}^{k}\{N_{n}^{X}(B_{i})=l_{i}\}, \frac{S_{n}}{\sqrt{n}}\leq y \right)\\
  &&= P\left(\bigcap_{i=1}^{k}\{N_{n}^{Z}(B_{i})=l_{i}\}, \frac{S_{n}}{\sqrt{n}}\leq y \right)\\
  &&= P\left(\bigcap_{i=1}^{k}\{N_{n}^{Z}(B_{i})=l_{i}\}\right)-
     P\left(\bigcap_{i=1}^{k}\{N_{n}^{Z}(B_{i})=l_{i}\}, \frac{S_{n}}{\sqrt{n}}> y \right) \\
  &&= P\left(\bigcap_{i=1}^{k}\{N_{n}^{Z}(B_{i})=l_{i}\}\right)-
     P\left(\bigcap_{i=1}^{k}\{N_{n}^{Y}(B_{i})=l_{i}\}, \frac{S_{n}}{\sqrt{n}}> y \right)  \\
  &&= P\left(\bigcap_{i=1}^{k}\{N_{n}^{Z}(B_{i})=l_{i}\}\right)-
     P\left(\bigcap_{i=1}^{k}\{N_{n}^{Y}(B_{i})=l_{i}\}\right)P\left(\frac{S_{n}}{\sqrt{n}}> y \right)
\end{eqnarray*}
where $l_{1},l_{2},\cdots,l_{k}$ are non-negative integer numbers.
Since $\{Y_{n}, n\geq1\}$  is an independent copy of $\{X_{n}, n\geq1\}$,
$N_{n}^{Y}$ converges in distribution to a Poisson process on (0,1] with parameter $\ln G^{-1}(x)$. Thus we have
$$P\left(\bigcap_{i=1}^{k}\{N_{n}^{Y}(B_{i})=l_{i}\}\right)\rightarrow\prod_{i=1}^{k}\frac{(m(B_{i})\log G^{-1}(x))^{l_{i}}}
{l_{i}!}e^{-m(B_{i})\log G^{-1}(x)},$$
as $n\rightarrow\infty$. By Lemma 2.9
$$P\left(\bigcap_{i=1}^{k}\{N_{n}^{Z}(B_{i})=l_{i}\}\right)\rightarrow\prod_{i=1}^{k}\frac{(m(B_{i})\log G^{-1}(x))^{l_{i}}}
{l_{i}!}e^{-m(B_{i})\log G^{-1}(x)},$$
as $n\rightarrow\infty$. Thus,
$$P\left(\bigcap_{i=1}^{k}\{N_{n}^{X}(B_{i})=l_{i}\}, \frac{S_{n}}{\sqrt{n}}\leq y \right)\rightarrow
\Phi(y)\prod_{i=1}^{k}\frac{(m(B_{i})\log G^{-1}(x))^{l_{i}}}{l_{i}!}e^{-m(B_{i})\log G^{-1}(x)},$$
as $n\rightarrow\infty$, which proves the asymptotic independence between $N_{n}^{X}$
and $S_{n}$. The proof of Theorem 2.1 is complete.

\textbf{Proof of Theorem 2.3:}
It is easy to see that
\begin{eqnarray*}
P\left(a_{n}(M_{n}^{(k)}-b_{n})\leq
x, \frac{S_{n}}{\sqrt{n}}\leq y\right)
                     &=&P(N_{n}^{X}((0,1])\leq k-1, \frac{S_{n}}{\sqrt{n}}\leq y)\\
                     &=&\sum_{l=0}^{k-1}P(N_{n}^{X}((0,1])= l, \frac{S_{n}}{\sqrt{n}}\leq y)
\end{eqnarray*}
and the result follows by Theorem 2.1.

\textbf{Proof of Theorem 2.4:} The proof of Theorem 2.4 is essentially the same as that of Theorem 2.1, and the details
are omitted here.

\textbf{Proof of Theorem 2.5:} By Theorem 2.4, we have as $n\rightarrow\infty$
\begin{eqnarray*}
&&P\left(a_{n}(M_{n}^{(1)}-b_{n})\leq x_{1},\ldots, a_{n}(M_{n}^{(k)}-b_{n})\leq x_{k}, \frac{S_{n}}{\sqrt{n}}\leq y\right)\\
&&=P\left(N_{n}^{(1)}((0,1])\leq 0,\ldots, N_{n}^{(k)}((0,1])\leq k-1, \frac{S_{n}}{\sqrt{n}}\leq y\right)\\
&&\rightarrow P\left(N^{(1)}((0,1])\leq 0,\ldots, N^{(k)}((0,1])\leq k-1\right)\Phi(y),
\end{eqnarray*}
which shows the asymptotic independence of the $k$ largest maxima and the partial sum.

\textbf{Proof of Corollary 2.6:}
By Theorem 2.4 again, we have
\begin{eqnarray*}
&&P\left(a_{n}(M_{n}^{(k)}-b_{n})\leq x_{1}, a_{n}(M_{n}^{(l)}-b_{n})\leq x_{2}, \frac{S_{n}}{\sqrt{n}}\leq y\right)\\
&&=P\left(N_{n}^{(k)}((0,1])\leq k-1, N_{n}^{(l)}((0,1])\leq l-1, \frac{S_{n}}{\sqrt{n}}\leq y\right)\\
&&\rightarrow P\left(N^{(k)}((0,1])\leq k-1, N^{(l)}((0,1])\leq l-1\right)\Phi(y),
\end{eqnarray*}
as $n\rightarrow\infty$. Note that
\begin{eqnarray*}
&&P\left(N^{(k)}((0,1])\leq k-1, N^{(l)}((0,1])\leq l-1\right)\\
&&=\sum_{i=0}^{k-1}\sum_{j=i}^{l-1}P\left(N^{(k)}((0,1])=i, N^{(l)}((0,1])=j\right)\\
&&=\sum_{i=0}^{k-1}\sum_{j=i}^{l-1}P\left(N^{(l)}((0,1])=j\right)P\left(N^{(k)}((0,1])=i\big|N^{(l)}((0,1])=j\right)\\
&&=\sum_{i=0}^{k-1}\sum_{j=i}^{l-1}\frac{\left(\log G^{-1}(x_{2})\right)^{j}}{j!}e^{-\log G^{-1}(x_{2})}\frac{j!}{i!(j-i)!}\left(\frac{\log G^{-1}(x_{1})}{\log G^{-1}(x_{2})}\right)^{i}
\left(1-\frac{\log G^{-1}(x_{1})}{\log G^{-1}(x_{2})}\right)^{j-i}\\
&&=\sum_{i=0}^{k-1}\sum_{j=i}^{l-1} G(x_{2})\frac{(\log G^{-1}(x_{1}))^{i}(\log G^{-1}(x_{2})-\log G^{-1}(x_{1}))^{j-i}}{i!(j-i)!},
\end{eqnarray*}
which completes the proof.

\section{The almost sure limit theorem}
 In this section, we extend Theorems 2.3 and 2.5 to the almost sure version.

\textbf{Theorem 3.1}. {\sl Let $\{X_{n}, n\geq1\}$  be a sequence of standard i.i.d. random variables with nondegenerate common distribution function $F$ and $E(X_{1})=0$, $E(X_{1}^{2})=1$. Suppose that there exists constants $a_{n}>0,b_{n}\in \mathbb{R}, n\ge1$ and a non-degenerate distribution $G(x)$ such that
(\ref{eq2.1})
 holds. Then for any $x,y\in \mathbb{R}$ and fixed $k\geq1$
\begin{equation}
\label{eq2.3}
 \lim_{n\rightarrow \infty}\frac{1}{\log
N}\sum^{N}_{n=1}\frac{1}{n}I\left(a_{n}(M_{n}^{(k)}-b_{n})\le
x,\frac{S_{n}}{\sqrt{n}}\leq y\right)=G(x)\sum_{i=0}^{k-1}\frac{(-\ln G(x))^{i}}{i!}\Phi(y)\ \ a.s.
\end{equation}
and for any $x_{1},\ldots,x_{k},y\in \mathbb{R}$ and fixed $k\geq1$
\begin{eqnarray}
\label{eq2.4}
&& \lim_{n\rightarrow \infty}\frac{1}{\log
N}\sum^{N}_{n=1}\frac{1}{n}I\left(a_{n}(M_{n}^{(1)}-b_{n})\leq x_{1},\ldots, a_{n}(M_{n}^{(k)}-b_{n})\leq x_{k},\frac{ S_{n}}{\sqrt{n}}\leq y\right)\nonumber\\
&&=H(x_{1},x_{2},\ldots,x_{k})\Phi(y)\ \ a.s.
\end{eqnarray}
}

Before giving the proof of Theorem 3.1, we state two lemmas which will be used in the proofs of our main results.
For $n-m\geq k$, let $M_{m,n}^{(k)}$  be the $k$-th maximum of $\{X_{m+1},\ldots,X_{n}\}$ and $S_{m,n}=\sum_{l=m+1}^{n}X_{l}$.
As usual, $a_{n}\ll b_{n}$ means $\limsup_{n\rightarrow\infty}|a_{n}/b_{n}|<+\infty$.

\textbf{Lemma 3.2}. {\sl Let $(\xi_{k})^{\infty}_{k=1}$ be a sequence
of uniformly bounded random variables, i.e., there exists some
$M\in (0,\infty)$ such that $|\xi_{k}|\leq M$  $a.s.$ for all $k\in
N$. If
\begin{eqnarray*}
Var\left(\sum^{N}_{n=1}\frac{1}{n}\xi_{n}\right) \ll
(\log N)^{2}(\log\log N)^{-(1+\varepsilon)}
\end{eqnarray*}
for some $\varepsilon>0$, then
\begin{eqnarray*}
\lim_{N\rightarrow\infty}\frac{1}{\log
N}\sum^{N}_{n=1}\frac{1}{n}(\xi_{n}-E\xi_{n})=0\ \ a.s.
\end{eqnarray*}
}

\textbf{Proof}. See Lemma $3.1$ of Cs\'{a}ki and
Gonchigdanzan (2002).

\textbf{Lemma 3.3}. {\sl Let $\{X_{n}, n\geq1\}$  be a sequence of standard i.i.d. random variables. We have for any fixed $k$ and $m\geq k$
$$P\left(M_{m}^{(k)}>M_{m,n}^{(k)}\right)\leq k\frac{m}{n}.$$
}

\textbf{Proof}. See Lemma 1 of Peng et al. (2009).

\textbf{Proof of Theorem 3.1}. Let
$$\eta_{n}=I\left(M_{n}^{(k)}\leq u_{n}, \frac{S_{n}}{\sqrt{n}}\leq y\right)-P\left(M_{n}^{(k)}\leq u_{n}, \frac{S_{n}}{\sqrt{n}}\leq y\right).$$
Notice that $(\eta_{n})_{n=1}^{\infty}$ is a sequence of bounded
random variables with $Var(\eta_{n})\leq 1$.
We first show
\begin{equation}
\label{eq3.7}
\lim_{N\rightarrow\infty}\frac{1}{\log N}\sum_{n=1}^{N}\frac{1}{n}\eta_{n}=0,\ \  a.s.
\end{equation}
Using Lemma 3.2, we only need to show
\begin{equation}
\label{eq3.8}
Var\left(\sum^{N}_{n=1}\frac{1}{n}\eta_{n}\right)\ll(\log N)^{2}(\log\log N)^{-(1+\varepsilon)}.
\end{equation}
We have,
\begin{eqnarray*}
Var\left(\sum^{N}_{n=1}\frac{1}{n}\eta_{n}\right)&=& E\left(\sum^{N}_{n=1}\frac{1}{n}\eta_{n}\right)^{2} \\
                                      &=& \sum^{N} _{n=1}\frac{E\eta_{n}^{2}}{n^{2}}
                                        + 2\sum_{1\leq m<n\leq N}\frac{E(\eta_{m}\eta_{n})}{nm}\\
                                      &=:& L_{N,1}+2L_{N,2}.
\end{eqnarray*}
Clearly
\begin{eqnarray*}
L_{N,1}= \sum^{N} _{n=1}\frac{1}{n^{2}}E\eta_{n}^{2}
    \leq  \sum^{N} _{n=1}\frac{1}{n^{2}}
    =  O(1).
\end{eqnarray*}
For $L_{N,2}$, for $n\geq m+k$, we have
\begin{eqnarray*}
\left|E(\eta_{m}\eta_{n})\right|&=&\left|Cov\left(I\left(M_{m}^{(k)}\leq u_{m}, \frac{S_{m}}{\sqrt{m}}\leq y\right),I\left(M_{n}^{(k)}\leq u_{n},\frac{S_{n}}{\sqrt{n}}\leq y\right)\right)\right|\\
&\leq& R_{1}+R_{2}+R_{3},
\end{eqnarray*}
where
\begin{eqnarray*}
R_{1}=\left|Cov\left(I\left(M_{m}^{(k)}\leq u_{m}, \frac{S_{m}}{\sqrt{m}}\leq y\right),\left[I\left(M_{n}^{(k)}\leq u_{n},\frac{S_{n}}{\sqrt{n}}\leq y\right)-I\left(M_{m,n}^{(k)}\leq u_{n},\frac{S_{n}}{\sqrt{n}}\leq y\right)\right]\right)\right|,
\end{eqnarray*}
\begin{eqnarray*}
R_{2}=\left|Cov\left(I\left(M_{m}^{(k)}\leq u_{m}, \frac{S_{m}}{\sqrt{m}}\leq y\right),\left[I\left(M_{m,n}^{(k)}\leq u_{n},\frac{S_{n}}{\sqrt{n}}\leq y\right)-I\left(M_{m,n}^{(k)}\leq u_{n},\frac{S_{m,n}}{\sqrt{n}}\leq y\right)\right]\right)\right|
\end{eqnarray*}
and
\begin{eqnarray*}
R_{3}=\left|Cov\left(I\left(M_{m}^{(k)}\leq u_{m}, \frac{S_{m}}{\sqrt{m}}\leq y\right),I\left(M_{m,n}^{(k)}\leq u_{n},\frac{S_{m,n}}{\sqrt{n}}\leq y\right)\right)\right|.
\end{eqnarray*}
For the first term $R_{1}$, by Lemma 3.3, we have
\begin{eqnarray*}
R_{1}&\ll& E\left|I\left(M_{n}^{(k)}\leq u_{n},\frac{S_{n}}{\sqrt{n}}\leq y\right)-I\left(M_{m,n}^{(k)}\leq u_{n},\frac{S_{n}}{\sqrt{n}}\leq y\right)\right|\\
&\ll& P\left(M_{m,n}^{(k)}\leq u_{n}\right)-P\left(M_{n}^{(k)}\leq u_{n}\right)\\
&\ll& P\left(M_{m}^{(k)}> M_{m,n}^{(k)}\right)\\
&\leq& k\frac{m}{n}.
\end{eqnarray*}
For the second term $R_{2}$, we have
\begin{eqnarray*}
R_{2}&\ll& E\left|I\left(M_{m,n}^{(k)}\leq u_{n},\frac{S_{n}}{\sqrt{n}}\leq y\right)-I\left(M_{m,n}^{(k)}\leq u_{n},\frac{S_{m,n}}{\sqrt{n}}\leq y\right)\right|\\
&\ll& E\left|\frac{S_{m}}{\sqrt{n}}\right|\\
&\leq& \frac{1}{\sqrt{n}}(E(S_{m}^{2}))^{1/2}\\
&=& \left(\frac{m}{n}\right)^{1/2}.
\end{eqnarray*}
By the independence of  $\{X_{n}, n\geq 1\}$, we have
$$R_{3}=0.$$
Then we conclude that
\begin{eqnarray*}
L_{N,2}&\ll&\sum_{1\leq m<n\leq N\atop n\geq m+k}\frac{1}{mn}\left(\frac{m}{n}+\left(\frac{m}{n}\right)^{1/2}\right)+
\sum_{1\leq m<n\leq N\atop n< m+k}\frac{1}{mn}\\
&\leq& \sum_{n=2}^{N}\frac{1}{n^{2}}\sum_{m=1}^{n-1}1 + \sum_{n=2}^{N}\frac{1}{n^{3/2}}\sum_{m=1}^{n-1}\frac{1}{m^{1/2}}+\sum_{m=1}^{N-1}\frac{1}{m}\sum_{n=m+1}^{m+k}\frac{1}{n}\\
 &\leq& \sum_{n=2}^{N}\frac{1}{n} + 2\sum_{n=2}^{N}\frac{1}{n}+\sum_{m=1}^{N-1}\frac{1}{m}\\
&\ll&\log N.
\end{eqnarray*}
Thus, (\ref{eq3.8}) holds. Note that Theorem 2.3 implies
\begin{equation}
\label{eq3.9}
\lim_{N\rightarrow\infty}\frac{1}{\log N}\sum_{n=1}^{N}\frac{1}{n}P\left(M_{n}^{(k)}\leq u_{n}, \frac{S_{n}}{\sqrt{n}}\leq y\right)=G(x)\sum_{i=0}^{k-1}\frac{(-\ln G(x))^{i}}{i!}\Phi(y)\ \ a.s.,
\end{equation}
and then the first assertion of Theorem 2.1 follows from (\ref{eq3.7}) and (\ref{eq3.9}).
The second assertion can be proved similarly.

{\bf Declarations}

{\bf Competing interests}

The authors declare no competing interests.

{\bf Author contributions}

Li was a major contributor in writing the manuscript; Tan provided some helpful discussions in writing the manuscript.
All authors read and approved the manuscript.



\begin{thebibliography}{100} \small

\bibitem{}Anderson, C. W., Turkman, K. F. The joint limiting distribution of sums and maxima of stationary
sequences. J. Appl. Probab., 1991a 28, 33-44.

\bibitem{}Anderson, C. W., Turkman, K. F. Sums and maxima in stationary sequences, J Appl Probab, 1991b,
 28, 715-716.

\bibitem{}Chow, T. L., Teugels, J. L. The sum and the maximum of i.i.d. random variables. In: Proceedings of the
Second Prague Symposium on Asymptotic Statistics (Hradec Kralove, 1978), North-Holland, New York,
1979, 81-92

\bibitem{} Dudzi\'{n}ski, M. An almost sure limit theorem for the maxima and sums of stationary Gaussian
sequences. Probability and Mathematical Statistics, 2003, 23, 139-152.

\bibitem{} Dudzi\'{n}ski, M. The almost sure central limit theorems in the joint version for the maxima and
sums of certain stationary Gaussian sequences. Statistics and Probability Letters, 2008, 78, 347-357.

\bibitem{}Ho, H. C., Hsing, T. On the asymptotic joint distribution of the sum and maximum of stationary normal
random variables. J. Appl. Probab., 1996, 33, 138-145.

\bibitem{} Ho, H.C., McCormick, W.P. Asymptotic distribution of sum and maximum for strongly dependent
Gaussian processes. J. Appl. Probab., 1999, 36, 1031-1044.

\bibitem{}Hsing, T. A note on the asymptotic independence of the sum and maximum of strongly mixing stationary
random variables. Ann. Probab., 1995, 23, 938-947.


\bibitem{}Hu, A., Peng, Z., Qi, Y. Joint behavior of point process of exceedances and partial sum from a
Gaussian sequence. Metrika, 2009, 70, 279-295.

\bibitem{}James, B., James, K., Qi, Y. Limit distribution of the sum and maximum from multivariate Gaussian
sequences. J. Multivariate Anal., 2007, 98, 517-532.

\bibitem{}Leadbetter, M. R., Lindgren, G., Rootz\'{e}n, H. Extremes and Related Properties of Random Sequences and
Processes, Springer, New York, 1983.


\bibitem{} McCormick, W.P., Qi, Y. Asymptotic distribution for the sum and the maximum of Gaussian processes.
J. Appl. Probab., 2000, 37, 958-971.


\bibitem{} Peng, Z. Joint asymptotic distributions of exceedances point process and partial sum of strong dependent
Gaussian sequences. Acta Math. Appl. Sin., 1998, 22, 362-367.

\bibitem{} Peng, Z., Li, J., Nadarajah, S.  Almost sure convergence of extreme
order statistics. Electronic Journal of Statistics, 2009, 3, 546-556.

\bibitem{} Peng, Z., Nadarajah, S. On the joint limiting distribution of sums and maxima of stationary normal
sequence. Theory Probab. Appl., 2002, 47, 817-820.


\bibitem{} Peng, Z., Tong, J., Weng, Z. Joint limit distributions of exceedances point processes and partial sums of gaussian vector sequence. Acta Mathematica Sinica. English Series, 2012, 28(8): 1647-1662.

\bibitem{} Peng, Z., Wang, L., Nadarajah, S. Almost sure central limit
theorem for partial sums and maxima. Mathematische Nachrichten, 2009, 282, 632-636.



\bibitem{}Tan, Z., Peng, Z.  Joint asymptotic distributions of exceedances point process and partial sum of  strong dependent
nonstationary  Gaussian sequences. Acta Math. Appl. Sin., 2011a, 34, 24-32.

\bibitem{}Tan, Z., Peng, Z. Joint asymptotic distribution of exceedances point process and partial sum of stationary Gaussian sequence. Appl. Math. J. Chinese Univ. 2011b, 26, 319-326.

\bibitem{}Tan, Z.,  Wang, Y., Almost sure central limit theorem for the maxima and sums of stationary Gaussian sequences,  Journal of the Korean Statistical Society, 2011, 40(3), 347-355.

\bibitem{} Tan, Z., Yang, Y.,  The maxima and sums of multivariate non-stationary Gaussian sequences, Appl. Math. J. Chinese Univ., 2015, 30, 197-209.

\bibitem{} Wu, Q. Improved results in almost sure central limit theorems for the maxima and partial sums for Gaussian sequences, Journal of Inequalities and Applications, 2015, article number: 109.

\bibitem{} Wu, Q., Jiang, Y.  Almost sure central limit theorem for self-normalized partial sums and maxima. Revista de la Real Academia de Ciencias Exactas, Fisicas y Naturales, Serie A. Matematicas, 2016, 110(2), 699-710.


\bibitem{} Zang, Q.P. A note on the almost sure central limit theorems for the maxima and sums. Journal of Inequalities and Applications, 2012, article number: 223.

\bibitem{} Zang, Q.P.,  Wang, Z. and Fu, K.  A note on almost sure central limit theorem in the joint version for the maxima and sums. Journal of Inequalities and Applications, 2010, article number: 234964.

\bibitem{} Zhao, S., Peng, Z., Wu, S. Almost sure convergence for the maximum and the sum of nonstationary guassian sequences. Journal of Inequalities and  Applications, 2010, article number: 856495

\end{thebibliography}
\end{document}